\documentclass{article}
\usepackage{amsmath, amsthm, amsfonts, amssymb}
\usepackage{hyperref}
\usepackage[latin1]{inputenc}

\makeatletter \@addtoreset{equation}{section}

\makeatletter

\def\beg   {\begin{theorem}}   \def\ee   {\end{theorem}}
\def\be   {\begin{equation}}   \def\ee   {\end{equation}}
\def\ba   {\begin{array}}      \def\ea   {\end{array}}
\def\bea  {\begin{eqnarray}}   \def\eea  {\end{eqnarray}}
\def\bean {\begin{eqnarray*}}  \def\eean {\end{eqnarray*}}

\newtheorem{lemma}{Lemma}[section]
\newtheorem{theorem} [lemma]{Theorem}

\newtheorem{remark}[lemma]{Remark}

\newtheorem{example}[lemma] {Example}

\newcommand{\mC}{\ensuremath{\mathbb{C}}}

\newcommand{\mN}{\ensuremath{\mathbb{N}}}

\begin{document}

\vspace{4cm}
\begin{center} \LARGE{\textbf{ Uniqueness of $p(f)$ and $P[f]$ }}
 \end{center}
 \vspace{1cm}
 \begin{center} \bf{Kuldeep Singh Charak$^{1 }$, \quad Banarsi Lal$^{2}$ }
\end{center}

\begin{center}
 Department of Mathematics, University of Jammu,
Jammu-180 006, INDIA.\\
{$^{1}$ E-mail: kscharak7@rediffmail.com }\\

{$^{2}$ E-mail: banarsiverma644@gmail.com }
\end{center}

\bigskip
\begin{abstract}
Let $f$ be a non-constant meromorphic function, $a (\not\equiv 0, \infty)$ be a meromorphic function satisfying $T(r,a) = o(T(r,f))$ as $r \rightarrow \infty$, and $p(z)$ be a polynomial of degree $n \geq 1$ with $p(0) = 0$. Let $P[f]$ be a non-constant differential polynomial of $f$. Under certain essential conditions, we prove the uniqueness of $p(f)$ and $P[f]$  when $p(f)$ and $P[f]$ share $a$ with weight $l \geq 0$. Our result generalizes the results due to Zang and Lu, Banerjee and Majumdar, Bhoosnurmath and Kabbur and answers a question of Zang and Lu.
\end{abstract}

\vspace{1cm}\noindent \textbf{Keywords: } Meromorphic functions, small functions, sharing of values, differential polynomials, Nevanlinna theory.

\vspace{0.5cm} \noindent\textbf{AMS subject classification: 30D35, 30D30}

\vspace{4cm}

\normalsize
\newpage

\section{Introduction}

   Let $f$ and $g$ be two non constant meromorphic functions and $k$ be a non-negative integer. For $a \in \mC \cup \{\infty \}$, we denote by $E_k(a, f)$ the set of all a-points of $f$, where an a-point of multiplicity $m$ is counted $m$ times if $m \leq k$ and $k +1$ times if $m > k$. If $E_k(a, f) = E_k(a, g)$, we say that $f$ and $g$ share the value $a$ with weight $k$.\\
 
 We write ``$f$ and $g$ share $(a, k)$" to mean that ``$f$ and $g$ share the value $a$ with weight $k$". Since $E_k(a, f) = E_k(a, g)$ implies $E_p(a, f) = E_p(a, g)$ for any integer $p (0 \leq p <k)$, clearly if $f$ and $g$ share $(a, k)$, then $f$ and $g$ share $(a, p)$, $0 \leq p <k$. Also we note that $f$ and $g$ share the value $a$ IM(ignoring multilicity) or CM(counting multiplicity) if and only if $f$ and $g$ share $(a, 0)$ or $(a, \infty)$, respectively.
 
\medskip

 A {\it differential polynomial} $P[f]$ of a non-constant meromorphic function $f$ is defined as 
 $$P[f]: = \sum^{m}_{i = 1}M_i[f],$$
where $M_i[f] = a_i.\prod^{k}_{j = 0}(f^{(j)})^{n_{ij}}$ with $n_{i0},\;n_{i1},\dotsc, n_{ik}$ as non-negative integers and $a_i (\not\equiv 0)$ are meromorphic functions satisfying $T(r,a_i) = o(T(r,f))$ as $r \rightarrow \infty$. The numbers $\overline d(P) = {\text max}_{1 \leq i \leq m}\sum^k_{j = 0}n_{ij}$ and $\underline d(P) = {\text min}_{1 \leq i \leq m}\sum^k_{j = 0}n_{ij}$ are respectively called the degree and lower degree of $P[f]$. If $\overline d(P) = \underline d(P) = d$ (say), then we say that  $P[f]$ is a homogeneous differential polynomial of degree $d$.
 
\medskip

For notational purpose, let $f$ and $g$ share 1 IM, and let $z_0$ be a zero of $f - 1$ with multiplicity $p$ and a zero of $g - 1$ with multiplicity $q$. We denote by $N^{1)}_E\left(r, 1/(f - 1)\right)$, the counting function of the zeros of $f - 1$ when $p = q = 1$.  By $\overline N^{(2}_E\left(r, 1/(f - 1)\right)$, we denote the counting function of the zeros of $f - 1$ when $p = q \geq 2 $ and by $\overline N_L\left(r, 1/(f - 1)\right)$, we denote the counting function of the zeros of $f - 1$ when $p > q \geq 1$, each point in these counting functions is counted only once; similarly, the terms $N^{1)}_E\left(r, 1/(g - 1)\right)$, $\overline N^{(2}_E\left(r, 1/(g - 1)\right)$ and $\overline N_L\left(r, 1/(g - 1)\right)$. Also, we denote by $\overline N_{f > k}\left(r, 1/(g - 1)\right)$, the reduced counting function of those zeros of $f - 1$ and $g - 1$ such that $p > q = k$, and similarly the term  $\overline N_{g > k}\left(r, 1/(f - 1)\right)$.\\

\medskip

 Inspired by a uniqueness result due to Mues and Steinmetz \cite{MS} : ``{\it If $f$ is a non-constant entire function sharing two distinct values ignoring multiplicity with $f'$, then $f \equiv f'$} ", the study of the uniqueness of $f$ and $f^{(k)}$, $f^n$ and $(f^m)^{(k)}$, $f$ and $P[f]$ is carried out by numerous authors. For example, Zang and Lu \cite{TL} proved :
 
\medskip

{\bf Theorem A.} \ {\it Let $k$, $n$ be the positive integers, $f$ be a non-constant meromorphic function, and $a (\not\equiv 0, \infty)$ be a meromorphic function satisfying $T(r,a) = o(T(r,f))$ as $r \rightarrow \infty$. If $f^n$  and $f^{(k)}$ share $a$ IM and 
\begin{equation*}
(2k + 6)\Theta(\infty, f) + 4 \Theta(0, f) + 2 \delta _{2 + k}(0,f) > 2k + 12 -n,
\end{equation*}
or $f^n$  and $f^{(k)}$ share $a$ CM and 
\begin{equation*} \label{2}
(k + 3)\Theta(\infty, f) + 2 \Theta(0, f) +  \delta _{2 + k}(0,f) > k + 6 -n,
\end{equation*}
then $f^n \equiv f^{(k)}$.}

\medskip 

In the same paper, T. Zhang and W. Lu asked the following question:

\medskip

{\bf Question 1:} {\it What will happen if $f^n$ and $(f^{(k)})^m$ share a meromorphic function $a (\not\equiv 0, \infty)$ satisfying $T(r,a) = o(T(r,f))$ as $r \rightarrow \infty$ ?} \\

 S.S.Bhoosnurmath  and Kabbur \cite{BK} proved:
 
 \medskip
 
{\bf Theorem B.} \ {\it Let $f$ be a non-constant meromorphic function and $a (\not\equiv 0, \infty)$ be a meromorphic function satisfying $T(r,a) = o(T(r,f))$ as $r \rightarrow \infty$. Let $P[f]$ be a non-constant differential polynomial of $f$.
If $f$  and $P[f]$ share $a$ IM and 
\begin{equation*}
(2Q + 6)\Theta(\infty, f) + (2 + 3 \underline d(P)) \delta(0,f) > 2Q + 2 \underline d(P) + \overline d(P) + 7,
\end{equation*}
or if $f$ and $P[f]$ share $a$ CM and 
\begin{equation*}
3\Theta(\infty, f) + (\underline d(P) + 1)\delta(0,f) > 4,
\end{equation*}
then $f \equiv P[f]$.}

\medskip

Banerjee and Majumder \cite{BM} considered the weighted sharing of $f^n$ and $(f^m)^{(k)}$ and proved the following result:

\medskip

{\bf Theorem C.} \ {\it Let $f$ be a non-constant meromorphic function, $k,n,m \in \mN$ and $l$ be a non negative integer. Suppose $a (\not\equiv 0, \infty)$ be a meromorphic function satisfying $T(r,a) = o(T(r,f))$ as $r \rightarrow \infty$ such that $f^n$  and $(f^m)^{(k)}$ share $(a, l)$. If $l \geq 2$ and 
\begin{equation*} 
(k + 3)\Theta(\infty, f) + (k + 4) \Theta(0, f) > 2k + 7 -n,
\end{equation*} 
or $l = 1$ and 
\begin{equation*} 
\left(k + \frac{7}{2}\right)\Theta(\infty, f) + \left(k + \frac{9}{2}\right)\Theta(0, f) > 2k + 8 -n,
\end{equation*} 
or $l = 0$ and 
\begin{equation*}
(2k + 6)\Theta(\infty, f) + (2k + 7) \Theta(0, f) > 4k + 13 -n,
\end{equation*}
then $f^n \equiv (f^m)^{(k)}$.}

\bigskip

Motivated by such uniqueness investigations, it is rational to think about the problem in more general setting: {\it Let $f$ be a non-constant meromorphic function, $P[f]$ be a non-cnstant differential polynomial of $f,$ $p(z)$ be a polynomial of degree $n \geq 1$ and $a (\not\equiv 0, \infty)$ be a meromorphic function satisfying $T(r,a) = o(T(r,f))$ as $r \rightarrow \infty$. If $p(f)$ and $P[f]$ share $(a,l), \ l \geq 0$, then is it true that $p(f) \equiv P[f]\;?$}\\

Generally this is not true, but under certain essential conditions, we prove the following result:

\begin{theorem}
Let $f$ be a non-constant meromorphic function, $a (\not\equiv 0, \infty)$ be a meromorphic function satisfying $T(r,a) = o(T(r,f))$ as $r \rightarrow \infty$, and $p(z)$ be a polynomial of degree $n \geq 1$ with $p(0) = 0$. Let $P[f]$ be a non-constant differential polynomial of $f$. Suppose $p(f)$  and $P[f]$ share $(a, l)$ with one of the following conditions: \\ 
$(i)$ $l \geq 2$ and 
\begin{equation} \label{1}
(Q + 3)\Theta(\infty, f) + 2n \Theta(0, p(f)) + \overline d(P) \delta(0,f) > Q + 3 + 2 \overline d(P) - \underline d(P) + n,
\end{equation} 
$(ii)$ $l = 1$ and 
\begin{equation} \label{2}
\left(Q + \frac{7}{2}\right)\Theta(\infty, f) + \frac{5n}{2}\Theta(0, p(f)) + \overline d(P)\delta(0,f) > Q + \frac{7}{2} + 2\overline d(P) - \underline d(P) + \frac{3n}{2},
\end{equation}
$(iii)$ $l = 0$ and 
\begin{equation} \label{3}
(2Q + 6)\Theta(\infty, f) + 4n \Theta(0, p(f)) + 2\overline d(P)\delta(0,f) > 2Q + 6 + 4 \overline d(P) - 2\underline d(P) + 3n.
\end{equation}
Then $p(f) \equiv P[f].$
\end{theorem}

\begin{example} Consider the function $f(z) = cos \alpha z + 1 - 1/{\alpha ^4}$, where $\alpha \neq 0,\pm1, \pm i$ and $p(z) = z$. Then $p(f)$ and $ P[f] \equiv f^{(iv)}$ share $(1, l)$, $l \geq 0$ and none of the inequalities (\ref{1}), (\ref {2}) and (\ref{3}) is satisfied, and $p(f) \neq P[f]$. Thus conditions in Theorem $1.1$ can not be removed.
\end{example}

\begin{remark} Theorem $1.1$ generalizes Theorem $A$, Theorem $B$, Theorem $C$ (and also generalizes Theorem $1.1$ and Theorem $1.2$ of \cite{BM}) and provides an answer to a question of Zhang and Lu \cite{TL}.
\end{remark}

The main tool of our investigations in this paper is Nevanlinna value distribution theory\cite{HAY}.

\section{Proof of the Main Result}
 We shall use the following results in the proof of our main result: 
 \begin{lemma} \label{a}\cite{BK}
 Let $f$ be a non-constant meromorphic function and $P[f]$ be a differential polynomial of $f$. Then
 \begin{equation} \label{u}
 m \left(r, \frac{P[f]}{f^{\overline d(P)}}\right) \leq (\overline d(P) - \underline d(P))m\left(r, \frac{1}{f}\right) + S(r,f),
 \end{equation}
 \begin{equation} \label{v}
 N \left(r, \frac{P[f]}{f^{\overline d(P)}}\right) \leq (\overline d(P) - \underline d(P))N\left(r, \frac{1}{f}\right) + Q \left[\overline N(r, f) + \overline N\left(r, \frac{1}{f}\right)\right] + S(r,f),
 \end{equation}
 \begin{equation} \label{w}
 N\left(r, \frac{1}{P[f]}\right) \leq Q \overline N(r,f) + (\overline d(P) - \underline d(P))m\left(r, \frac{1}{f}\right)+ N\left(r, \frac{1}{f^{\overline d(P)}}\right) + S(r,f),
 \end{equation}
  where $Q = {\text max}_{1 \leq i \leq m} \{n_{i0}+ n_{i1} + 2n_{i2} + ... +kn_{ik}\}$.
\end{lemma}

\begin{lemma} \label{b} \cite{AB}
Let $f$ and $g$ be two non-constant meromorphic functions.\\
{\it (i)} If $f$ and $g$ share $(1, 0)$, then
\begin{align}\label{y}
\overline N_L\left(r, \frac{1}{f-1}\right) \leq \overline N\left(r, \frac{1}{f}\right) + \overline N(r,f) + S(r),
\end{align}
where $S(r) = o(T(r))$ as $r \rightarrow \infty$ with $T(r) = {\text max}\{T(r,f);T(r,g)\}$.\\
{\it (ii)} If $f$ and $g$ share $(1, 1)$, then 
\begin{align}\label{z}
2\overline N_L\left(r, \frac{1}{f-1}\right) + 2 \overline N_L\left(r, \frac{1}{g-1}\right)
& + \overline N^{(2}_E\left(r, \frac{1}{f-1}\right) - \overline N_{f > 2}\left(r, \frac{1}{g-1}\right)\notag\\
& \leq N\left(r, \frac{1}{g - 1}\right) - \overline N\left(r, \frac{1}{g - 1}\right).
\end{align}
\end{lemma}

\bigskip 

{\bf Proof of Theorem 1.1:}
Let $F = p(f)/a$ and $G = P[f]/a$. Then
\begin{equation} \label{4}
F - 1 = \frac{p(f) - a}{a}\; \text{and}\;G - 1 = \frac{P[f] - a}{a}.
\end{equation}
Since $p(f)$  and $P[f]$ share $(a, l)$, it follows that $F$ and $G$ share $(1, l)$ except at the zeros and poles of $a$. Also note that 
$$\overline N(r, F) = \overline N(r, f) + S(r,f) \;\text{and}\;\overline N(r, G) = \overline N(r, f) + S(r,f).$$

Define
\begin{equation} \label{5}
\psi = \left(\frac{F''}{F'} - \frac{2F'}{F - 1}\right) - \left(\frac{G''}{G'} - \frac{2G'}{G - 1}\right).
\end{equation}
{\bf Claim:} $\psi \equiv 0$.\\

Suppose on the contrary that $\psi \not\equiv 0$. Then from (\ref{5}), we have
 $$m(r, \psi) = S(r, f).$$
 By the Second fundamental theorem of Nevanlinna, we have
\begin{align} \label{9}
T(r, F) + T(r, G) & \leq  2\overline N(r, f) + \overline N\left(r, \frac{1}{F}\right) + \overline N\left(r, \frac{1}{F - 1}\right) + \overline N\left(r, \frac{1}{G}\right)\notag\\
&+ \overline N\left(r, \frac{1}{G - 1}\right) -  N_0\left(r, \frac{1}{F'}\right) - N_0\left(r, \frac{1}{G'}\right) + S(r, f),
\end{align}
where  $ N_0(r, 1/F')$ denotes the counting function of the zeros of $F'$ which are not the zeros of $F(F - 1)$ and $ N_0(r, 1/G')$ denotes the counting function of the zeros of $G'$ which are not the zeros of $G(G - 1)$.

\medskip

 {\bf Case 1.} When $l \geq 1$.
 
 \medskip
 
  Then from (\ref{5}), we have,
 \begin{align*}
N^{1)}_E\left(r, \frac{1}{F - 1}\right) & \leq  N \left(r,\frac{1}{\psi} \right) + S(r,f) \notag\\
& \leq  T(r, \psi) + S(r,f) \notag\\
& =  N(r, \psi) + S(r,f) \notag\\
& \leq  \overline N(r, F) + \overline N_{(2}\left(r, \frac{1}{F}\right) + \overline N_{(2}\left(r, \frac{1}{G}\right) + \overline N_L\left(r, \frac{1}{F - 1}\right)\notag\\
&+ \overline N_L\left(r, \frac{1}{G - 1}\right) +  N_0\left(r, \frac{1}{F'}\right) + N_0\left(r, \frac{1}{G'}\right) + S(r, f).
\end{align*}
and so
\begin{align}\label{10}
\overline N\left(r, \frac{1}{F - 1}\right)+\overline N\left(r, \frac{1}{G - 1}\right)& = N^{1)}_E\left(r, \frac{1}{F - 1}\right)+\overline N^{(2}_E\left(r, \frac{1}{F - 1}\right) + \overline N_L\left(r, \frac{1}{F - 1}\right)\notag\\
&+ \overline N_L\left(r, \frac{1}{G - 1}\right) + \overline N\left(r, \frac{1}{G - 1}\right) + S(r,f) \notag\\
&\leq \overline N(r, f) + \overline N_{(2}\left(r, \frac{1}{F}\right) + \overline N_{(2}\left(r, \frac{1}{G}\right) + 2\overline N_L\left(r, \frac{1}{F - 1}\right)\notag\\
&+ 2 \overline N_L\left(r, \frac{1}{G - 1}\right) + \overline N^{(2}_E\left(r, \frac{1}{F - 1}\right) +  \overline N\left(r, \frac{1}{G - 1}\right)\notag\\ 
&+ N_0\left(r, \frac{1}{F'}\right) + N_0\left(r, \frac{1}{G'}\right) + S(r, f).
\end{align}
{\bf Subcase 1.1:} When $l = 1$.\\
In this case, we have
\begin{equation}\label{12}
\overline N_L\left(r, \frac{1}{F - 1}\right) \leq \frac{1}{2}N\left(r, \frac{1}{F'}|F \neq 0\right) \leq \frac{1}{2}\overline N(r, F) + \frac{1}{2}\overline N\left(r, \frac{1}{F}\right),
\end{equation}
where $N\left(r, \frac{1}{F'}|F \neq 0\right)$ denotes the zeros of $F'$, that are not the zeros of $F$.

\medskip

 From (\ref{z}) and (\ref{12}), we have
\begin{align} \label{13}
2\overline N_L\left(r, \frac{1}{F - 1}\right) + 2 \overline N_L\left(r, \frac{1}{G - 1}\right)& + \overline N^{(2}_E\left(r, \frac{1}{F - 1}\right) + \overline N\left(r, \frac{1}{G - 1}\right)\notag\\
&\leq N\left(r, \frac{1}{G - 1}\right) + \overline N_L\left(r, \frac{1}{F - 1}\right) + S(r, f)\notag \\
&\leq N\left(r, \frac{1}{G - 1}\right) + \frac{1}{2} \overline N(r, F) + \frac{1}{2}\overline N\left(r, \frac{1}{F}\right) + S(r,f)\notag \\
&\leq N\left(r, \frac{1}{G - 1}\right) + \frac{1}{2} \overline N(r, f) + \frac{1}{2}\overline N\left(r, \frac{1}{p(f)}\right) + S(r,f).
\end{align}
Thus, from (\ref{10}) and (\ref{13}), we have
\begin{align} \label{14}
\overline N\left(r, \frac{1}{F - 1}\right)+\overline N\left(r, \frac{1}{G - 1}\right) & \leq  \overline N(r, f) + \overline N_{(2}\left(r, \frac{1}{F}\right) + \overline N_{(2}\left(r, \frac{1}{G}\right)\notag\\
&+ \frac{1}{2} \overline N(r, f) + \frac{1}{2}\overline N\left(r, \frac{1}{p(f)}\right) + N\left(r, \frac{1}{G - 1}\right)\notag \\
&+ N_0\left(r, \frac{1}{F'}\right) + N_0\left(r, \frac{1}{G'}\right) + S(r, f)\notag\\
& \leq  \overline N(r, f) + \overline N_{(2}\left(r, \frac{1}{F}\right) + \overline N_{(2}\left(r, \frac{1}{G}\right)\notag\\
&+ \frac{1}{2} \overline N(r, f) + \frac{1}{2}\overline N\left(r, \frac{1}{p(f)}\right) + T(r, G)\notag \\
&+ N_0\left(r, \frac{1}{F'}\right) + N_0\left(r, \frac{1}{G'}\right) + S(r, f).
\end{align}
From (\ref{w}), (\ref{9}) and (\ref{14}), we obtain
\begin{align*}
T(r,F) & \leq 3\overline N(r, f) + \overline N\left(r, \frac{1}{F}\right) + \overline N_{(2}\left(r, \frac{1}{F}\right) +  \overline N\left(r, \frac{1}{G}\right)+
\overline N_{(2}\left(r, \frac{1}{G}\right)\notag\\
& + \frac{1}{2} \overline N(r, f) + \frac{1}{2}\overline N\left(r, \frac{1}{p(f)}\right)+ S(r,f)\notag\\
& \leq \frac{7}{2}\overline N(r, f) + 2\overline N\left(r, \frac{1}{F}\right) + N\left(r, \frac{1}{G}\right) + \frac{1}{2}\overline N\left(r, \frac{1}{p(f)}\right) + S(r,f)\notag\\
& \leq \frac{7}{2}\overline N(r, f) + \frac{5}{2}\overline N\left(r, \frac{1}{p(f)}\right) + N\left(r, \frac{1}{P[f]}\right) + S(r,f)\notag\\
&\leq \left(Q + \frac{7}{2}\right)\overline N(r, f) + \frac{5}{2}\overline N\left(r, \frac{1}{p(f)}\right) + (\overline d(P) - \underline d(P)) T(r,f) +\overline d(P)N\left(r, \frac{1}{f}\right) + S(r,f)\notag\\
&\leq \left[\left(Q + \frac{7}{2}\right)\{1 - \Theta (\infty, f)\} + \frac{5n}{2} \{1 - \Theta (0, p(f))\} +\overline d(P) \{1 - \delta (0, f)\}\right] T(r, f)\notag\\
& + (\overline d(P)-\underline d(P)) T(r,f) + S(r, f).
\end{align*}
That is,
\begin{align*}
n T(r,f) & =  T(r,F) + S(r,f)\\
&\leq \left[\left(Q + \frac{7}{2}\right)\{1 - \Theta (\infty, f)\} + \frac{5n}{2} \{1 - \Theta (0, p(f))\} +\overline d(P) \{1 - \delta (0, f)\}\right] T(r, f)\notag\\
& + (\overline d(P)-\underline d(P)) T(r,f) + S(r, f).
\end{align*}
Thus 
$$\left[\{\left(Q + \frac{7}{2}\right)\Theta(\infty, f) + \frac{5n}{2}\Theta(0, p(f)) + \overline d(P)\delta(0,f)\} - \{Q + \frac{7}{2} + 2\overline d(P) - \underline d(P) + \frac{3n}{2}\}\right] T(r, f) \leq S(r,f).$$
That is,
$$\left(Q + \frac{7}{2}\right)\Theta(\infty, f) + \frac{5n}{2}\Theta(0, p(f)) + \overline d(P)\delta(0,f) \leq Q + \frac{7}{2} + 2\overline d(P) - \underline d(P) + \frac{3n}{2},$$
which violates (\ref{2}).

\medskip

{\bf Subcase 1.2:} When $l \geq 2$.\\
In this case, we have  
$$2\overline N_L\left(r, \frac{1}{F - 1}\right) + 2 \overline N_L\left(r, \frac{1}{G - 1}\right) + \overline N^{(2}_E\left(r, \frac{1}{F - 1}\right) + \overline N\left(r, \frac{1}{G - 1}\right) \leq N\left(r, \frac{1}{G - 1}\right) + S(r, f).$$
Thus from (\ref{10}), we obtain
\begin{align} \label{11}
\overline N\left(r, \frac{1}{F - 1}\right)+\overline N\left(r, \frac{1}{G - 1}\right)& \leq \overline N(r, f) + \overline N_{(2}\left(r, \frac{1}{F}\right) + \overline N_{(2}\left(r, \frac{1}{G}\right) +  N\left(r, \frac{1}{G - 1}\right)\notag\\
&+ N_0\left(r, \frac{1}{F'}\right) + N_0\left(r, \frac{1}{G'}\right) + S(r, f)\notag\\
& \leq \overline N(r, f) + \overline N_{(2}\left(r, \frac{1}{F}\right) + \overline N_{(2}\left(r, \frac{1}{G}\right) +  T(r, G)\notag\\
&+ N_0\left(r, \frac{1}{F'}\right) + N_0\left(r, \frac{1}{G'}\right) + S(r, f).
 \end{align}
Now from (\ref{w}), (\ref{9}) and (\ref{11}), we obtain
\begin{align*}
T(r,F) &\leq  3\overline N(r, f) + \overline N\left(r, \frac{1}{F}\right) + \overline N_{(2}\left(r, \frac{1}{F}\right) + \overline N\left(r, \frac{1}{G}\right) + \overline N_{(2}\left(r, \frac{1}{G}\right) + S(r,f) \notag\\
& \leq 3\overline N(r, f) + 2\overline N\left(r, \frac{1}{F}\right) + N\left(r, \frac{1}{G}\right) + S(r,f) \notag\\
& \leq 3\overline N(r, f) + 2\overline N\left(r, \frac{1}{p(f)}\right) + N\left(r, \frac{1}{P[f]}\right) + S(r,f)\notag\\
&\leq (Q + 3)\overline N(r, f) + 2\overline N\left(r, \frac{1}{p(f)}\right) + (\overline d(P)-\underline d(P)) T(r,f) +\overline d(P)N\left(r, \frac{1}{f}\right) + S(r,f)\notag\\
&\leq [(Q + 3)\{1 - \Theta (\infty, f)\} + 2n \{1 - \Theta (0, p(f))\} +\overline d(P) \{1 - \delta (0, f)\}] T(r, f)\notag\\
& + (\overline d(P)-\underline d(P)) T(r,f) + S(r, f).
\end{align*}
 That is,
 \begin{align*}
nT(r,f) & =  T(r, F) + S(r,f)\\
&\leq [(Q + 3)\{1 - \Theta (\infty, f)\} + 2n \{1 - \Theta (0, p(f))\} +\overline d(P) \{1 - \delta (0, f)\}] T(r, f)\notag\\
& + (\overline d(P)-\underline d(P)) T(r,f) + S(r, f).
\end{align*}
Thus
$$[\{(Q + 3)\Theta(\infty, f) + 2n \Theta(0, p(f)) + \overline d(P) \delta(0,f)\} - \{(Q + 3 + 2 \overline d(P) - \underline d(P) + n\}] T(r,f) \leq S(r,f).$$
That is, 
$$(Q + 3)\Theta(\infty, f) + 2n \Theta(0, p(f)) + \overline d(P) \delta(0,f) \leq Q + 3 + 2 \overline d(P) - \underline d(P) + n,$$
which violates (\ref{1}).\\

{\bf Case 2.} When $l = 0$.

\medskip

 Then, we have
\begin{equation*} 
N^{1)}_E\left(r, \frac{1}{F - 1}\right) = N^{1)}_E\left(r, \frac{1}{G - 1}\right) + S(r, f),\; \overline N^{(2}_E\left(r, \frac{1}{F - 1}\right) = \overline N^{(2}_E\left(r, \frac{1}{G - 1}\right) + S(r, f),
\end{equation*} 
and also from (\ref{5}), we have
\begin{align}\label{15}
\overline N\left(r, \frac{1}{F - 1}\right)+\overline N\left(r, \frac{1}{G - 1}\right)&\leq N^{1)}_E\left(r, \frac{1}{F - 1}\right)+\overline N^{(2}_E\left(r, \frac{1}{F - 1}\right) + \overline N_L\left(r, \frac{1}{F - 1}\right)\notag\\
&+ \overline N_L\left(r, \frac{1}{G - 1}\right) + \overline N\left(r, \frac{1}{G - 1}\right) + S(r,f) \notag\\
& \leq  N^{1)}_E\left(r, \frac{1}{F - 1}\right) + \overline N_L\left(r, \frac{1}{F - 1}\right) + N\left(r, \frac{1}{G - 1}\right) + S(r,f) \notag\\
& \leq  \overline N(r, F) + \overline N_{(2}\left(r, \frac{1}{F}\right) + \overline N_{(2}\left(r, \frac{1}{G}\right) + 2\overline N_L\left(r, \frac{1}{F - 1}\right)\notag\\
&+ \overline N_L\left(r, \frac{1}{G - 1}\right) + N\left(r, \frac{1}{G - 1}\right) + N_0\left(r, \frac{1}{F'}\right)\notag\\
&+ N_0\left(r, \frac{1}{G'}\right) + S(r, f).
\end{align}
 From (\ref{w}),(\ref{y}),(\ref{9}) and (\ref{15}), we obtain
\begin{align*}
T(r,F) &\leq  3\overline N(r, f) + \overline N\left(r, \frac{1}{F}\right) + \overline N_{(2}\left(r, \frac{1}{F}\right) + \overline N\left(r, \frac{1}{G}\right) + \overline N_{(2}\left(r, \frac{1}{G}\right)\notag\\
&+ 2\overline N_L\left(r, \frac{1}{F - 1}\right) + \overline N_L\left(r, \frac{1}{G - 1}\right) + S(r,f) \notag\\
& \leq  3\overline N(r, f) + 2\overline N\left(r, \frac{1}{F}\right) + N\left(r, \frac{1}{G}\right) + 2\overline N\left(r, \frac{1}{F}\right)\notag\\
&+ 2 \overline N(r, F) + \overline N\left(r, \frac{1}{G}\right) +  \overline N(r, G) + S(r,f) \notag\\
& \leq  6\overline N(r, f) + 4\overline N\left(r, \frac{1}{F}\right) + 2N\left(r, \frac{1}{G}\right) + S(r,f) \notag\\
& \leq 6\overline N(r, f) + 4\overline N\left(r, \frac{1}{p(f)}\right) + 2N\left(r, \frac{1}{P[f]}\right) + S(r,f)\notag \\
&\leq (2Q + 6)\overline N(r, f) + 4\overline N\left(r, \frac{1}{p(f)}\right) + 2(\overline d(P) - \underline d(P)) T(r,f) + 2 \overline d(P)N\left(r, \frac{1}{f}\right) + S(r,f)\notag\\
&\leq [(2Q + 6)\{1 - \Theta (\infty, f)\} + 4n \{1 - \Theta (0, p(f))\} + 2\overline d(P) \{1 - \delta (0, f)\}] T(r, f)\notag\\
& + 2(\overline d(P)-\underline d(P)) T(r,f) + S(r, f).
\end{align*}
That is,
\begin{align*}
n T(r,f) & = T(r,F) + S(r,f)\\
&\leq [(2Q + 6)\{1 - \Theta (\infty, f)\} + 4n \{1 - \Theta (0, p(f))\} + 2\overline d(P) \{1 - \delta (0, f)\}] T(r, f)\notag\\
& + 2(\overline d(P)-\underline d(P)) T(r,f) + S(r, f).
\end{align*}
Thus 
$$[\{(2Q + 6)\Theta(\infty, f) + 4n \Theta(0, p(f)) + 2\overline d(P)\delta(0,f)\} - \{2Q + 6 + 4 \overline d(P) - 2\underline d(P) + 3n\}] T(r, f) \leq S(r,f).$$
That is,
$$(2Q + 6)\Theta(\infty, f) + 4n \Theta(0, p(f)) + 2\overline d(P)\delta(0,f) \leq 2Q + 6 + 4 \overline d(P) - 2\underline d(P) + 3n,$$
which violates (\ref{3}).

\medskip

This proves the claim and thus $\psi \equiv 0$. So (\ref{5}) implies that
$$\frac{F''}{F'} - \frac{2F'}{F - 1} = \frac{G''}{G'} - \frac{2G'}{G - 1},$$
and so we obtain
\begin{equation}\label{6}
\frac{1}{F - 1} = \frac{C}{G - 1} + D,
\end{equation} 
where $C \neq 0$ and $D$ are constants.\\

Here, the following three cases can arise:

\medskip

{\bf Case$(i):$} When $D \neq 0,\;-1$. Rewriting (\ref{6}) as
$$\frac{G - 1}{C} = \frac{F - 1}{D + 1 - DF},$$
we have
$$\overline N(r, G) = \overline N\left(r,\frac{1}{F - (D + 1)/D}\right).$$

\medskip

In this subcase, the Second fundamental theorem of Nevanlinna yields
\begin{align*}
nT(r,f) & =  T(r, F) + S(r,f)\\
& \leq  \overline N(r,F) + \overline N\left(r,\frac{1}{F}\right) + \overline N\left(r,\frac{1}{F - (D + 1)/D}\right) + S(r,f)\\
& \leq  \overline N(r,F) + \overline N\left(r,\frac{1}{F}\right) + \overline N(r, G) + S(r,f)\\
& \leq  2\overline N(r,f) + \overline N\left(r,\frac{1}{p(f)}\right) + S(r,f)\\
& = [2 \{1 - \Theta(\infty, f)\} + n \{1 - \Theta(0, p(f))\}] T(r, f) + S(r, f).
\end{align*}
Thus 
$$[\{2\Theta(\infty, f) +  n \Theta(0, p(f))\} - 2] T(r, f) \leq S(r, f).$$
That is,
$$2\Theta(\infty, f) +  n \Theta(0, p(f)) \leq 2,$$
which contradicts (\ref{1}),(\ref{2}) and (\ref{3}).

\medskip 

{\bf Case$(ii):$} When $D = 0$. Then from (\ref{6}), we have
\begin{equation}\label{7}
G = CF - (C - 1).
\end{equation}
So if $C \neq 1$, then
$$\overline N\left(r, \frac{1}{G}\right) = \overline N\left(r,\frac{1}{F - (C - 1)/C}\right).$$

\medskip

Now the Second fundamental theorem of Nevanlinna and (\ref{w}) gives
\begin{align*}
nT(r,f) & =  T(r, F) + S(r,f)\\
& \leq  \overline N(r,F) + \overline N\left(r,\frac{1}{F}\right) + \overline N\left(r,\frac{1}{F - (C - 1)/C}\right) + S(r,f)\\
& \leq  \overline N(r,F) + \overline N\left(r,\frac{1}{F}\right) + \overline N\left(r, \frac{1}{G}\right) + S(r,f)\\
& \leq  \overline N(r,f) + \overline N\left(r,\frac{1}{p(f)}\right) + \overline N\left(r, \frac{1}{P[f]}\right) + S(r,f)\\
& \leq  \overline N(r,f) + \overline N\left(r, \frac{1}{p(f)}\right) + Q \overline N(r, f) + (\overline d(P) - \underline d(P)) m\left(r,\frac{1}{f}\right)\\
&+ N\left(r, \frac{1}{f^{\overline d(P)}}\right) + S(r,f)\\
& \leq  (Q + 1)\overline N(r,f) + \overline N\left(r, \frac{1}{p(f)}\right) + (\overline d(P) - \underline d(P)) T(r,f)\\
&+ \overline d(P)N\left(r, \frac{1}{f}\right) + S(r,f)\\
&\leq  [(Q + 1) \{1 - \Theta(\infty, f)\} + n \{1 - \Theta(0, p(f))\} + \overline d(P) \{1 - \delta(0,f)\}] T(r, f)\\
&+ (\overline d(P) - \underline d(P)) T(r,f) + S(r, f).
\end{align*}
Thus 
$$[\{(Q + 1)\Theta(\infty, f) + n \Theta(0, p(f)) + \overline d(P)\delta(0,f)\} - \{Q + 1 + 2 \overline d(P) - \underline d(P)\}] T(r, f) \leq S(r, f).$$
That is,
$$(Q + 1)\Theta(\infty, f) + n \Theta(0, p(f)) + \overline d(P)\delta(0,f) \leq Q + 1 + 2 \overline d(P) - \underline d(P),$$
which contradicts (\ref{1}),(\ref{2}) and (\ref{3}).

\medskip

Thus, $C = 1$ and so in this case from (\ref{7}), we obtain $F \equiv G$ and so
$$p(f) \equiv P[f].$$
{\bf Case$(iii):$} When $D = -1$. Then from (\ref{6}) we have
\begin{equation} \label{8}
\frac{1}{F - 1} = \frac{C}{G - 1} - 1.
\end{equation}
So if $C \neq -1$, then
$$\overline N\left(r, \frac{1}{G}\right) = \overline N\left(r,\frac{1}{F - C/(C + 1)}\right),$$
and as in the Subacase $(ii)$, we find that
\begin{align*}
nT(r,f) & \leq  (Q + 1)\overline N(r,f) + \overline N\left(r, \frac{1}{p(f)}\right) + (\overline d(P) - \underline d(P)) T(r,f)\\
&+ \overline d(P)N\left(r, \frac{1}{f}\right) + S(r,f).
\end{align*}
Thus 
$$[\{(Q + 1)\Theta(\infty, f) + n \Theta(0, p(f)) + \overline d(P)\delta(0,f)\} - \{Q + 1 + 2 \overline d(P) - \underline d(P)\}]T(r, f) \leq S(r, f).$$
That is,
$$(Q + 1)\Theta(\infty, f) + n \Theta(0, p(f)) + \overline d(P)\delta(0,f) \leq Q + 1 + 2 \overline d(P) - \underline d(P),$$
which contradicts (\ref{1}),(\ref{2}) and (\ref{3}).
\medskip

Thus, $C = -1$ and so in this case from (\ref{8}), we obtain $FG \equiv 1$ and so $p(f)P[f] = a^2.$ Thus, in this case $\overline N(r,f) + \overline N\left(r, 1/f\right) = S(r,f)$.

\medskip

 Now, by using (\ref{u}) and (\ref{v}), we have
\begin{align*}
(n + \overline d(P))T(r,f) &\leq  T\left(r, \frac{a^2}{f^{n + \overline d(P)}}\right) + S(r,f)\\
& \leq  T\left(r, \left[1 + \frac{a_{n - 1}}{f} + --- + \frac{a_1}{f^{n-1}} \right ].\frac{P[f]}{f^{\overline d(P)}}\right) + S(r,f)\\
& \leq (n-1)T(r, f) + T\left(r, \frac{P[f]}{f^{\overline d(P)}} \right) + S(r, f)\\
&= (n-1)T(r, f) +  m\left(r, \frac{P[f]}{f^{\overline d(P)}}\right) + N\left(r, \frac{P[f]}{f^{\overline d(P)}}\right) + S(r,f)\\
&\leq (n-1)T(r, f) + (\overline d(P) - \underline d(P))m\left(r, \frac{1}{f}\right) + (\overline d(P) - \underline d(P))N\left(r, \frac{1}{f}\right)\\
&+ Q\left[\overline N(r,f) + \overline N\left(r, \frac{1}{f}\right)\right] + S(r,f)\\
& \leq (n-1)T(r, f) + (\overline d(P) - \underline d(P)) T(r, f) + S(r, f).
\end{align*}
Thus
$$(1 + \underline d(P)) T(r, f) \leq S(r, f),$$
which is a contradiction.

~~~~~~~~~~~~~~~~~~~~~~~~~~~~~~~~~~~~~~~~~~~~~~~~~~~~~~~~~~~~~~~~~~~~~~~~~~~~~~~~~~~~~~~~~~~~~~~~~~$\Box$

\bibliographystyle{amsplain}

\end{document}